\documentstyle[12pt]{article}
\date{}
\def\direcciones{\small \begin{tabular}{cc}
Esteban Andruchow &Gustavo Corach y Demetrio Stojanoff\\
& \\
Instituto de Ciencias&Instituto Argentino de Matem\'atica\\
Univ. Nac de Gral. Sarmiento&Saavedra 15 3er piso\\
Roca 850, (1663) San Miguel   & (1083) Buenos Aires, Argentina\\
Pcia. de Buenos Aires, Argentina&gcorach@mate.dm.uba.ar\\
eandruch@percanta.ungs.edu.ar&demetrio@mate.dm.uba.ar \\
\end{tabular}}

\newtheorem{fed}{Definition}[section]

\newtheorem{lem}[fed]{Lemma}

\newtheorem{pro}[fed]{Proposition}
\newtheorem{rem}[fed]{Remark}

\newtheorem{num}[fed]{}
\def\dem{{\it Proof.\ }\rm}

\def\zC{\hbox{\rm C\hskip -5.4pt\vrule
height 8.0pt width 0.4pt depth -0pt\hskip 4.5pt}}

\def\bull{\vrule height 1.0ex width .4ex depth -.1ex }
\def\inc{\subseteq}
\def\QED{\bull}

\def\ben{\begin{enumerate}}
\def\een{\end{enumerate}}
\def\noi{\noindent}
\def\csta{C$^*$-algebra}

\def\A1{{\cal A}}
\def\B1{{\cal B}}
\def\qa{Q}

\def\qsa{Q\sb S}
\def\qpa{Q\sb p}
\def\inv{^{-1}}
\def\*a{\#\sb a}
\def\pa{P}
\def\pau{P_1}
\def\pb{P(\B1 )}
\def\paa{P\sb a }
\def\op{{\cal O}\sb p}
\def\m2a{M\sb 2(\A1 )}
\def\aut{{\cal S}}
\def\auta{{\cal S} \sb a }
\def\fp{\varphi \sb p}
\def\vfi{\varphi }
\def\fa{\varphi \sb a}
\def\fq{\varphi \sb q}

\def\proj{\Omega}
\def\fia{\varphi\sb \alpha}
\def\U{\cal U\sb {\cal A}}
\def\H{{\cal H}}

\def\za{{\cal Z} (\A1)}
\def\ce{conditional expectation}
\def\la{\lambda}
\def\eps{\varepsilon}
\title{{\sc GEOMETRY OF OBLIQUE PROJECTIONS
\footnote{Research partially supported by CONICET, 
ANPCYT and UBACYT (Argentina)}}}

\bigskip

\author{ E. Andruchow, G. Corach and D. Stojanoff}

\begin{document}
 \maketitle 
\vskip1truecm
\begin{abstract}{
Let $\A1$ be a unital \csta . Denote by $\pa$ 
the space of selfadjoint projections of $\A1$. 
We study the relationship between $\pa$  and the spaces of
projections $ \paa $ determined by the different involutions 
$\*a$ induced by positive invertible elements 
$a \in \A1$. The maps $\fp : \pa \to \paa $
sending $p$ to the unique $q \in \paa$ with the same range as $p$ 
and $\proj _a : \paa \to \pa$ sending $q$ to the unitary 
part of the polar decomposition of the symmetry $2q-1$
are shown to be diffeomorphisms. 
We characterize the pairs
of idempotents $q, r \in \A1$ with $\|q-r\|<1$ 
such that there exists a positive 
element $a\in \A1$ verifying that $q, r \in \paa$. 
In this case $q$ and $r$ can be
joined by an unique short geodesic along the space of idempotents 
$\qa$ of $\A1$.
}\end{abstract}

\newdimen\normalbaselineskip
\normalbaselineskip=16pt
\normalbaselines
\vskip 1truecm
\section{Introduction.} Let   $\H$ be a Hilbert space with scalar product $<,>$. For every 
bounded positive invertible operator $a : \H \to \H$ consider the scalar 
product $<,>_a $ given by
$$
<\xi , \eta>\sb a = <a\xi , \eta> \ , \quad \xi , \ \eta \in \H . 
$$
It is clear that $<,>_a$ induces a norm equivalent to 
the norm induced by $<,>$. With respect to the scalar
product   $<,>_a$, the adjoint of a bounded linear operator
$x : \H \to \H$ is
$$
x^{\*a} = a \inv x^* a.
$$
Thus, $x$ is $a$-selfadjoint if and only if
$$
ax = x^*a.
$$

Given a closed subspace $S$ of $\H$, denote by $p = P_S$ the orthogonal 
projection from $\H$ onto $S$ and, for any positive
operator $a$, denote by $\fp (a) $ the unique $a$-selfadjoint projection
with range $S$. In a recent paper, Z. Pasternak-Winiarski \cite{[PW1]} proves the
analyticity of the map $a \mapsto \fp(a)$ and calculates its Taylor
expansion. This study is relevant for understanding reproducing kernels of
Hilbert spaces of holomorphic $L^2$ sections of complex vector bundles and
the way they change when the measures and hermitian structures are deformed
(see \cite{[PW2]}, \cite{[PW3]}). This type of deformations appears in a natural way when 
studying quantization of systems where the phase space is a K\"ahler 
manifold (Odzijewicz \cite{[O1]}, \cite{[O2]}).

In this paper we pose Pasternak-Winiarski's problem in the \csta \ setting
and use the knowledge of the differential geometry of idempotents,
projections and positive invertible elements in order to get more general
results in a shorter way.

More precisely, let $\A1$ be a unital \csta , $G = G(\A1 ) $ the 
group of invertible elements of $\A1$, ${\cal U} = \U$ 
the unitary group of $\A1$, 
$ G^+ = \{ a \in G : a^* = a , \ a \ge 0 \}$ 
the space of positive invertible elements of $\A1$, 
$$
\qa = Q(\A1 )= \{ q \in \A1 : q^2 = q\} \quad \hbox{ and } \quad 
\pa = P(\A1 ) = \{ p \in \qa : p = p ^*  \},
$$
the spaces of idempotents and projections of $\A1$.
The nonselfadjoint elements of $\qa$ will be called oblique 
projections. It is well known that $\qa$ is a closed analytic 
submanifold of $\A1$, $\pa$ is a closed real analytic 
submanifold of $\qa$ and $G^+$ is an open submanifold of 
$$
\aut = \aut (\A1 ) = \{ b \in \A1 : b^{*} = b \},
$$
which is a closed real subspace of $\A1$ (see \cite{[PR1]}, 
\cite{[CPR1]} or \cite{[CPR3]} for
details).

We define a fibration
$$
\vfi : \pa \times G^+ \to \qa 
$$
which coincides, when $\A1 = L(\H)$, with the map $(p,a) \mapsto \fp(a)$, 
the unique $a$-selfadjoint projection
with the same range as $p$. This allows us to study the analyticity of 
Pasternak-Winiarski's map in both variables $p$, $a$. The rich geometry
of $\qa$, $\pa$ and $G^+$ give an amount of information which may be
useful in the problems that motivated \cite{[PW1]}. 

Along this note we use the fact that  
every $p \in \qa$ induces a representation 
$\alpha \sb p$ of elements of $\A1$ by $2 \times 2$ matrices given by
$$
\alpha \sb p (a) = 
\left( \begin{array}{cc} pap &pa(1-p) \\ (1-p)ap & (1-p)a(1-p)
\end{array} \right). 
$$
Under this homomorphism $p$ can be identified with 
$$
\left( \begin{array}{cc} 1\sb {p\A1 p}&0 \\ 0&0 \end{array} \right) 
= \left( \begin{array}{cc} 1 &0 \\ 0&0 \end{array} \right) . 
$$
and  all idempotents $q$ with the same range of $p$ have the form 
$$
q = \left( \begin{array}{cc} 1&x \\ 0&0 \end{array} \right). 
$$
for some $x \in p\A1 (1-p)$. This trivial remark shortens many proofs in a
drastic way and the analyticity of some maps (for example
$\vfi : \pa \times G^+ \to \qa $) follows immediately.

The contents of the paper are the following. Section 2 contains some 
preliminary material including the matrix representations mentioned above
and the description of the adjoint operation induced by each positive
invertible (element or operator) $a$. 

In section 3 we  study the map $\fp = \vfi (p, .) : G^+ \to \qa$, which is
Pasternak-Winiarski's map when $\A1$ is $L(\H)$ and $p$ is the orthogonal 
projection $P_S$ onto a closed subspace $S \inc \H$. For $a \in G^+$, let
$\paa = \paa (\A1 )$ denote the set of all $\*a$-selfadjoint projections. 
This is a subset of $\qa$ and section 4 starts a study of the relationship 
between $\pa = \pau$ and $\paa$ and the way they are located in $\qa$. 
In particular we show that $\fia = 
\vfi (., a): \pa \to \paa $ is a diffeomorphism and compute its tangent 
map. Another interest\-ing map is the following: for $q \in \paa$, $\eps = 
2q-1$ is a reflection, i.e. $\eps ^2 = 1$, which admits in 
$\A1$ a polar decomposition $\eps = \la \rho$, with $\la \in G^+$ and 
$\rho $ a unitary element of $\A1$. It is easy to see that $\rho = 
\rho^* = \rho \inv $ so that $ p = \frac12(\rho+1) \in \pa$. In section 5
we prove that the map $\proj _ a: \paa \to \pa $ given by $\proj _a (q) = 
p$ is a diffeomorphism and study the movement of $\pa$ given by the
composition $\proj _ a \circ \fia : \pa \to \pa $. We also
characterize the orbit of $p$ by these movements, i.e.
$$
\op := \{ r \in \pa : \proj \sb a \circ \fa (p) = r \ \hbox{ for some } \ 
a \in G^+ \} . 
$$
In recent years several papers have appeared which study length of curves 
in $\pa$ and $\qa$ (see \cite{[PR2]}, \cite{[Br]}, \cite{[Ph]}, \cite{[CPR1]}, 
\cite{[ACS]} for example). It is
known that $\pa$ and the fibres of $\proj :\qa\to \pa$ are geodesically 
complete and their geodesics are short curves (for convenient Finsler
metrics -see \cite{[CPR1]}).  For a fixed $p\in \pa$, let us call horizontal
(resp. vertical) those directions around $p$ which produce geodesics 
along $\pa$ (resp. along the fiber $\proj \inv (p)$). 
In section 6 we  show that there exist short geodesics in many other 
directions (not only the horizontal
and the vertical ones). 

This paper, which started from a close 
examination of Pasternak-Winiar\-ski's
work, is part of the program of understanding the structure of the space of
idempotent operators. For a sample of the vast bibliography on the subject 
the reader is referred to the papers by Afriat \cite{[Af]}, Kovarik \cite{[Ko]}, Zem\'anek 
\cite{[Ze]}, Porta-Recht \cite{[PR1]}, 
Gerisch \cite{[Ge]}, Corach \cite{[C]} and the references therein. 
Applications of oblique projections to complex, harmonic and functional 
analysis and statistics can be found in the papers by Kerzman and 
Stein \cite{[KS1]}, \cite{[KS2]}, Ptak \cite{[Pt]}, Coifman and Murray \cite{[CM]} and Mizel and Rao 
\cite{[MR]}, among others. 

\section{Preliminary results.}
Let   $\H$ be a Hilbert space, $\A1 \subset L( \H)$  a unital \csta \, ,   
$G = G(\A1) $ the group of invertible elements and  
$\U$ the unitary group of $\A1$, 

If $S$ is a closed subspace of $\H$ and $q$ is a bounded linear projection
onto $S$, then
\begin{equation}\label{2.1}
p = qq^*(1-(q-q^*)^2)^{-1} 
\end{equation}
is the unique sefadjoint projection onto $S$. Note that, by this formula,
$p \in \A1$ when $q \in \A1$. Several different formulas are known for $p$
(see \cite{[Ge]}, p. 294); perhaps the simplest one is the so-called 
Kerzman-Stein formula
\begin{equation}\label{2.2}
p = q(1+ q-q^*)^{-1} 
\end{equation}
(see \cite{[KS1]}, \cite{[KS2]} or \cite{[CM]}). However, for the present purposes, 
(\ref{2.1}) is more convenient.
\medskip
We denote by 
\begin{equation}\label{2.3}
\qa = Q(\A1) = \{ q \in \A1 : q^2 = q\} \; \hbox{and} \;  
\pa = P(\A1) = \{ p \in \A1 : p = p ^* = p^2  \}
\end{equation}
the spaces of idempotents and projections of $\A1$.
Given a fixed closed subspace $S$ of $\H$, we denote by 
\begin{equation}\label{qs}
Q_S= Q\sb S(\A1) = \{ q \in Q(\A1) : q(\H) = S\}
\end{equation}
the space of idempotents of $\A1$ with range $S$. Note that, by (\ref{2.1}), 
$Q\sb S$ is not empty if and only if the projection 
$p = p\sb S$ onto $S$ belongs to $\A1$. We shall make this assumption.

It is easy to see that two idempotents $q,r\in\qa$ have the 
same range  if and only if $qr=r$ and $rq=q$. 
Therefore the space $Q\sb S $ of (\ref{qs}) can be characterized as
$$
Q\sb S  =  \qpa  = \{q\in\qa :qp=p ,\  pq=q  \}.
$$
In what follows, we shall adopt this notation $\qpa$, emphasizing the role
of $p$, rather than $S$. This will enable us to simplify many computations. 
Moreover this operator algebraic viewpoint allows one to get the results 
below independently of the representation of $\A1$.

\bigskip
\noi Recall some facts about matrix representations. 
Every $p \in \qa$ induces a representation 
$\alpha \sb p$ of elements of $\A1$ by $2 \times 2$ matrices given by
\begin{equation}\label{2.5}
\alpha \sb p (a) = 
\left( \begin{array}{cc} pap &pa(1-p) \\ (1-p)ap & (1-p)a(1-p)
\end{array} \right). 
\end{equation}
If $p \in \pa$ the representation preserves the involution $^*$.
\noi 
For simplicity we shall identify $a$ with $\alpha \sb p (a)$ and $\A1$ with 
its image by $\alpha \sb p $.  Observe that, with this convention, 
\begin{equation}\label{2.6}
p = \left( \begin{array}{cc} 1\sb {p\A1 p}&0 \\ 0&0 \end{array} \right) 
= \left( \begin{array}{cc} 1 &0 \\ 0&0 \end{array} \right) . 
\end{equation}
Moreover, $q \in \qsa = \qpa$ if and only if there exists $x \in p\A1 (1-p) $ such that
\begin{equation}\label{2.7}
q = \left( \begin{array}{cc} 1&x \\ 0&0 \end{array} \right). 
\end{equation}
Indeed, let $q = 
\left( \begin{array}{cc} a&b \\ c&d \end{array} \right) \in \qpa$.
Then 
$$
\left( \begin{array}{cc} 1&0 \\ 0&0 \end{array} \right) = 
p = qp = 
 \left( \begin{array}{cc} a&b \\ c&d \end{array} \right) 
\left( \begin{array}{cc} 1&0 \\ 0&0 \end{array} \right)= 
\left( \begin{array}{cc} a&0 \\ c&0 \end{array} \right), 
$$
then  $a=1 $ and $c=0$. On the other hand
$$
q= pq = \left( \begin{array}{cc} 1&b \\ 0&0 \end{array} \right),
$$
then $d =0$ and $b $ can be  anything. We summarize this information in the 
following:

\bigskip
\begin{pro}\label{2.8}
The space $\qpa$ can be identified with
$p\A1(1-p)$ by means of the affine map 
\begin{equation}\label{2.9}
\qpa \to p\A1(1-p) \  , \quad q \mapsto q-p
\end{equation}
\end{pro}
\smallskip
\noi Proof. \rm  Clearly, the affine map defined in  (\ref{2.9})  is  
injective. By (\ref{2.7}) it is well defined and onto.  \QED 

\medskip
\noi In the Hilbert space $\H$, every scalar product which is equivalent to 
the original $< , >$ is determined by a unique positive invertible 
operator $a \in L(\H)$ by means of 
\begin{equation}\label{2.10}
<\xi , \eta>\sb a = <a\xi , \eta> \ , \quad \xi , \ \eta \in \H . 
\end{equation}
For this scalar product the adjoint $x^{\*a}$ of 
$x \in L(\H)$ is easily seen 
to be
\begin{equation}\label{2.11}
x^{\*a} = a\inv x^* a 
\end{equation}
where $*$ denotes the adjoint operation for the original scalar product.
Operators which are selfadjoint for some $\*a$ have been considered by 
Lax \cite{[L]} and Dieu\-donn\'e \cite{[D]}. A geometrical study of families of 
C$^*$-involutions has been done by Porta and Recht \cite{[PR3]}. 

Denote by $G^+ = G^+(\A1)$ the set of all positive invertible elements of 
$\A1$. Every $a \in G^+$ induces as in (\ref{2.11}) a continuous involution 
${\*a}$ on $\A1$ by means of $x^{\*a} = a\inv x^* a$, for $x\in \A1$. 
$\A1$ is a \csta \ with the involution ${\*a}$ and the corresponding norm
$\|x\|_a = \|a^{1/2}xa^{-1/2}\|$ for $x \in \A1$. 
The mapping $x \mapsto a^{-1/2}x a^{1/2}$ is an isometric isomorphism
of $(\A1, \| \ \|, *)$ onto $(\A1, \| \ \|_a, \*a)$.
In this setting, $\A1$  can be also represented by the inclusion 
map in $L(\H, <,>\sb a)$. 

Note that the map $a \to <,>\sb a \mapsto \*a$ is not one 
to one, since (\ref{2.11}) 
says that if $ a \in \zC . I$ then $\*a=*$. If we regard  
this map in $G^+$ with values in the set of 
involutions of $\A1$,  then two elements $a, b \in G^+$ 
with $a = b z$ for $z$ in the center of $\A1$, 
\begin{equation}\label{2.12}
\za = \{z\in \A1 : zc=cz , \quad \hbox{ for all } \quad c\in \A1\},  
\end{equation}
produce the same involution $\*a$.

\medskip
\begin{num}\label{2.13} \rm
Recall the properties of  the conditional expectation 
induced by a fix\-ed projection $p\in \pa$. Note that the set 
$\A1\sb p$ of elements of $\A1$ which commute 
with $p$ is the C$^*$-subalgebra  
of $\A1$ of diagonal matrices  in terms of the representation (\ref{2.5}). 
We denote by $E\sb p : \A1 \to \A1\sb p\subset \A1 $ the \ce \ defined by 
compressing to the diagonal:
$$
E\sb p (a) = pap + (1-p)a(1-p) =
 \left( \begin{array}{cc} pap &0 \\ 0&(1-p)a(1-p) \end{array} \right) 
\ , \quad a\in \A1 . 
$$
This expectation has the following well known properties (\cite{[St]}, 
Chapter 2): for all
$a \in \A1$, 
\ben
\item  $E\sb p(bac) = bE\sb p(a) c $ for all $b,c \in \A1\sb p$.
\item $E\sb p (a^*) = E\sb p(a)^*.$
\item If $b\le a$  then   $E\sb p(b)\le E\sb p(a)$. 
In particular $E\sb p(G^+) \subset G^+$.
\item $\|E\sb p(a)\| \le \|a\|$.
\item If  $0\le a $, then $2E\sb p(a) \ge a$. 
\een
\end{num}

\bigskip
\section{Idempotents with the same range.}
\medskip

The main purpose of this section is to describe, for a fixed $p\in\pa$,  
the map which sends 
each $a \in G^+$ into the unique $q \in \qpa$ which is $\*a$-selfadjoint. 
This problem was posed and solved  by \cite{[PW1]} 
when $\A1=L(\H)$. Here we use $2\times 2$ matrix arguments to give very
short proofs of the results of \cite{[PW1]}. Moreover we  generalize
these results and apply them to understand some aspects of 
the geometry of the space $\qa$.

Let us fix the notations: For each $a\in G^+$ denote by 
\begin{equation}\label{3.1}
\auta = \auta (\A1) = \{ b \in \A1 : b^{\*a} = b \},  
\end{equation}
the set of $\*a$-selfadjoint elements of $\A1$.
\medskip
\begin{fed}\label{3.2}\rm
Let $\A1$ be a \csta \ and $p \in \pa$ a fixed 
projection of $\A1$. We consider the map
$$
\fp : G^+ \to \qpa  \quad \hbox{ given by } 
\fp (a) = \hbox { the unique } q \in \qpa \cap \auta \ , \quad a\in G^+ .
$$
Note that existence and uniqueness of such $q$ follow from (\ref{2.1})
applied to the \csta \ $\A1$ with the star $\*a$.
\end{fed}
\medskip
\bigskip
\begin{pro}\label{3.3}
Let $\A1$ be a \csta \ and $p \in \pa$. 
Then, for all $a\in G^+(\A1)$,
\begin{equation}\label{3.4}
\fp (a) = p E\sb p(a)\inv a , 
\end{equation}
where $E\sb p$ is the \ce \ defined in (\ref{2.13}). 
In particular, 
$$
\|\fp(a)\|\le 2\ \|a\|\ \|a\inv\|.
$$
\end{pro}
\smallskip
\noi Proof. \rm  Suppose that, in matrix form, we have
$$
a = \left( \begin{array}{cc} a\sb 1 &a\sb 2  \\ a\sb 2^*& a\sb 3 \end{array} \right) \quad  
\hbox{ and then } \quad E\sb p(a) = 
\left( \begin{array}{cc} a\sb 1 &0  \\ 0& a\sb 3 \end{array} \right).
$$
Since $\fp(a) \in \qpa$, by (\ref{2.7}) there exists 
$x \in p\A1(1-p)$ such that $\fp(a) =
p+x$. On the other hand, by (\ref{2.11}), 
$p+x \in \auta $ if and only if $ a\inv (p+x)^* a = p+x $, i.e.  
$ (p+x^*)a  = a(p+x)$. In matrix form,
$$
(p+x^*)a = \left( \begin{array}{cc} 1 &0  \\ x^*&0 \end{array} \right)
\left( \begin{array}{cc} a\sb 1 &a\sb 2  \\ a\sb 2^*& a\sb 3 \end{array} \right) = 
\left( \begin{array}{cc} a\sb 1 &a\sb 2  \\ x^*a\sb 1&x^* a\sb 2 \end{array} \right) \quad 
\hbox{ and}
$$
$$
a(p+x) = \left( \begin{array}{cc} a\sb 1 &a\sb 2  \\ a\sb 2^*& a\sb 3 \end{array} \right) 
\left( \begin{array}{cc} 1 &x  \\ 0&0 \end{array} \right) =
\left( \begin{array}{cc} a\sb 1 &a\sb 1x  \\ a\sb 2^*& a\sb 2^* x \end{array} \right).
$$
Then $(p+x^*)a = a(p+x)$ if and only if $a\sb 2 =a\sb 1 x$. Note that 
$a\in G^+(\A1) $ implies that $a\sb 1 \in G^+(p\A1 p)$. Then  
\begin{equation}\label{3.5}
\fp (a) = 
\left( \begin{array}{cc} 1 &a\sb 1\inv a\sb 2  \\ 0&0 \end{array} \right) , 
\end{equation}
and now formula (\ref{3.4}) can be proved by easy computations.
Finally, since $2E\sb p(a)\ge a$, we deduce that $E\sb p(a)\inv \le 2a\inv $ and 
the inequality $\|\fp(a)\|\le 2\|a\|\|a\inv\|$ follows easily. \QED
\medskip
\begin{rem}\label{3.6}\rm
There is a  way to describe $\fp$ in terms of
(\ref{2.2}) with the star $\*a$. 
In this sense we obtain, for $p \in \pa$ and $a\in G^+$, 
$$
\fp (a ) = p (1 + p - a\inv p a ) \inv 
= p(a+ap-pa)\inv a .
$$
Clearly $a+ap-pa = E\sb p(a) +2 a\sb 2^* $ and one  obtains 
(\ref{3.4}), since $p(a+ap-pa)\inv = pE\sb p(a)\inv $.  
However it seems difficult to obtain bounds for $\|\fp(a)\|$ 
by using this approach. 
\end{rem}
\bigskip
\begin{num} \rm
Consider the space $G^+$ as an open set of 
$\aut =\aut (\A1) ={\cal S}\sb 1(\A1) $, the closed real subspace  of 
selfadjoint elements of $\A1$. Then the map $\fp: G^+ \to \A1$ 
is real analytic. Indeed, if $h \in \aut $ and $\|h\|<1$, 
then 
\begin{equation}\label{3.7}
\fp (1+h) = p (1+E\sb p(h))\inv (1+h) =
p \sum\sb {n=0}^\infty (-1)^{n} E\sb p(h)^n (1+h), 
\end{equation}
and this formula is clearly real analytic near $1$. More computations 
starting from (\ref{3.7}) give the more explicit formula
\begin{equation}\label{3.8}
\fp (1+h) = p + \sum\sb {n=1}^\infty (-1)^{n-1} (ph)^n(1-p), 
\end{equation}
again for all $h \in \aut $ with $\|h\|<1$. These computations are 
very similar to those appearing in the proof of Theorem 5.1 of \cite{[PW1]}. We 
include them for the sake of completeness. By (\ref{3.7}),
$$
\begin{array}{rl}
\fp(1+h) & =  \ \sum\sb {n=0}^\infty (-1)^{n} (php)^n(p+ph) \\
& \\         
& =  \ \sum\sb {n=0}^\infty (-1)^{n}(php)^n + 
\sum\sb {n=0}^\infty (-1)^{n} (ph)^{n+1} \\
&\\    
     & =  \ p + \sum\sb {n=1}^\infty (-1)^{n} (ph)^n p  + 
\sum\sb {n=1}^\infty (-1)^{n-1} (ph)^{n} \\
&\\   
      & =  \ p + \sum\sb {n=1}^\infty (-1)^{n-1} (ph)^n (1-p). \end{array}
$$
As a consequence (see also Theorem 3.1 of \cite{[PW1]}) 
the tangent map $(T \fp )\sb 1 : \aut \to \A1 $ is given by 
\begin{equation}\label{3.9}
(T \fp )\sb 1 (X) = pX(1-p) \quad \hbox { for } \quad X\in \aut  .
\end{equation}
Actually, by (2.8) $\qpa$ is an affine manifold parallel to the closed 
subspace $p\A1 (1-p)$ which can be also regarded as its  ``tangent'' space. 
In this sense $(T \fp )\sb 1 $ is just the natural compression 
of $\aut $ onto $p\A1 (1-p)$. 

Note that formulas (\ref{3.7}) and (\ref{3.8}) do not depend on the selected star in 
$\A1$. Using this fact, formula (\ref{3.8}) can be generalized to a power 
series around each $a \in G^+$ by using (\ref{3.8}) with the star 
$\*a$ at $q= \fp (a)$.  Indeed, note that for every $b \in G^+$, 
$<,>\sb b= (<,>_a)_{a\inv b}$ is induced from $<,>\sb a$ by $a\inv b$, which 
is $a$-positive. If $h \in \aut$ and $\|h\|<\|a\inv \|\inv$, then 
$a+h \in G^+$, $\|a\inv h \|\sb a = \|a^{-1/2} ha^{-1/2} \| \le 
\|h\| \|a\inv \|  <1 $ and 
\begin{equation}\label{3.10}
\fp ( a+ h) = \fq ( 1 + a\inv h ) = 
 q + \sum\sb {n=1}^\infty (-1)^{n-1} (qa\inv h)^n(1-q) ,
\end{equation}
showing the real analyticity  of $\fp$ in $G^+$ and also giving the way to 
compute the tangent map $(T \fp)\sb a $ at every $a\in G^+$.
 
Formulas (\ref{3.9}),  (\ref{3.10}) and their consequence, the 
real analyticity 
of $\fp$ for $\A1 = L(\H)$, are the main results of \cite{[PW1]}. 
Here we generalize these results to an arbitrary \csta \ $\A1$.
In the following section, we shall explore 
some of their interesting geometrical interpretations and applications.
\end{num}
\bigskip
\section{Differential geometry of \qa .}
\medskip
The space $\qa$ of all idempotents of a \csta \ (or, more generally, of a
Banach algebra) has a rich topological and geometrical structure, 
studied for example in \cite{[MR]}, \cite{[Ze]}, \cite{[Ge]}, \cite{[PR1]}, 
\cite{[CPR1]} and \cite{[CPR2]}.

We recall some facts on the structure of $\qa$ as a closed submanifold of 
$\A1$. The reader is referred to \cite{[CPR1]} and \cite{[CPR2]} for details. 
The tangent space of $\qa$ at $q$ is naturally identified to 
\begin{equation}\label{TQ}
\begin{array}{rl}
\{ X \in \A1 : qX+Xq = X \} & = \{X\in \A1 : qXq = (1-q)X(1-q) =0 \} \\
& \\
& = q\A1 (1-q) \oplus (1-q) \A1 q . 
\end{array}
\end{equation}
In terms of the matrix representation induced by $q$, 
\begin{equation}\label{4.1}
T(\qa)\sb q = \{ \left( \begin{array}{cc} 0&x \\ y&0 \end{array} \right) \in \A1 \} 
\end{equation}
The set $\pa $ is a real submanifold of $\qa$. The tangent space $(T\pa )_p $
at $p \in \pa$ is 
$$
\{ X \in \A1 : pX+Xp = X , \ X^* =X \} ,
$$
which in terms of the matrix representation induced by $p$ is
\begin{equation}\label{TP}
T(\pa ) \sb p = \{
\left( \begin{array}{cc} 0&x^* \\ x&0 \end{array} \right) \in \A1\} = 
T(\qa)\sb p \cap \aut . 
\end{equation}
The space $\qa$  (resp. $\pa$) is a discrete union of homogeneous 
spaces of $G$ (resp. $\U$) by means of the natural action 
\begin{equation}\label{accion}
G\times \qa \to \qa \quad \hbox {given by } \quad 
(g,q) \mapsto gqg\inv 
\end{equation}
(resp. $\U \times \pa \to \pa$, $ (u,p) \mapsto upu^*$).

There is a natural connection on $\qa$ (resp. $\pa$) which induces in 
the tangent bundle $T\qa$ (resp. $T\pa$) a linear connection. 
The geodesics of this connection, i.e. the curves $\gamma$ such that 
the covariant derivative of $\dot \gamma$ vanishes, can be computed. 
For $X\in (T\qa)_p$  (resp. $(T\pa)_p$), the unique geodesic 
$\gamma $ with $\gamma (0) = p $ and $\dot \gamma (0) = X$ is 
given by
$$
\gamma (t) = e^{tX'} p e^{-tX'} ,
$$
where $X' = [X,p] = Xp-pX$. 
Thus, the exponential map $\exp \sb p: T(\qa)\sb p \to \qa $ is given by 
\begin{equation}\label{4.2}
\exp \sb p(X) = e^{ X'}pe^{-X'} \ ,\quad \quad \hbox{ for } \quad 
X \in T(\qa)\sb p .
\end{equation}
\bigskip

\begin{pro}\label{4.3}
The inverse of  the  affine bijective map 
$$
\Gamma: \qpa \to p\A1(1-p) \  , \quad \Gamma(q) = q-p . 
$$
of (\ref{2.9}) is  the restriction of the exponential map at $p$ to 
the closed subspace $p\A1(1-p) \subset  T(\qa)\sb p$. 
That is, for $x \in p\A1(1-p)$,  $ \exp\sb p (x) = p+x \in \qpa$.
\end{pro}
\smallskip
\noi Proof. \rm  Let $x \in p\A1(1-p)$. Then 
$$
\begin{array}{rl}
\exp\sb p (x) & =  
\exp \sb p \left( \begin{array}{cc} 0&x \\ 0&0 \end{array} \right)  \\
          &   \\
           & =   
\exp \left( \begin{array}{cc} 0&-x \\ 0&0 \end{array} \right) \ p \ 
\exp \left( \begin{array}{cc} 0&x \\ 0&0 \end{array} \right)     \quad 
          \hbox{ by (\ref{4.2}) }\\
&      \\           
& =   \left( \begin{array}{cc} 1&-x \\ 0&1 \end{array} \right)  \ p \ 
\left( \begin{array}{cc} 1&x \\ 0&1 \end{array} \right)  \\
&      \\           
& =  \left( \begin{array}{cc} 1&x \\ 0&0 \end{array} \right)  = p+x .  \quad 
\hbox { \QED } \end{array} 
$$
\begin{rem} \rm The map $\fp$ of \ref{3.2} can be also described using 
Proposition \ref{4.3}.
In fact, consider the real analytic map 
$$
u_p : G^+ \to G(\A1 ) \quad \hbox{ given by } \quad 
u_p(a) = \exp( - p E_p(a)\inv a (1-p) ) ,  \ a \in G^+ .
$$
Then, by \ref{4.3},   $\fp(a) = u_p(a) p u_p(a)\inv $. This is an 
explicit formula of an invertible  element which conjugates $p$ 
with $\fp(a)$. This can be a useful tool for lifting curves of 
idempotents to curves of invertible elements of $\A1$.
\end{rem}
\bigskip
\noi Now we consider the map $\fp$ by letting $p$ vary in $\pa$:  
\begin{equation}\label{4.4}
\vfi : \pa \times G^+ \to \qa \quad \hbox { given by } \quad
\vfi (p,a) = \fp (a) = p E\sb p(a)\inv a  , 
\end{equation}
for $p \in \pa , \ a \in G^+$. Consider also the map $\phi : \qa \to \pa $ 
given by (\ref{2.1}): 
\begin{equation}\label{4.5}
\phi (q) = qq^*(1-(q-q^*)^2)^{-1} , \quad \hbox{ for } \quad q \in \qa . 
\end{equation}
This map $\phi$ assigns 
to any $q \in \qa$ the unique $p \in \pa$ with the same range as $q$.

\bigskip
\begin{pro}\label{4.6}
The map $\vfi : \pa \times G^+ \to \qa$ is a $C^\infty$ fibration. 
For $q\in \qa $, let $ p = \phi (q)$ and $x = q-p \in p\A1 (1-p)$. 
Then the fibre of $q$ is 
\begin{equation}\label{4.7}
\vfi \inv (q) = 
\{ (p ,\left( \begin{array}{cc} a\sb 1 &a\sb 1 x  \\ 
x^* a\sb 1& a\sb 3 \end{array} \right) ) : 
0< a\sb 1  \ \hbox{ and }  
\ x^* a\sb 1 x < a\sb 3    \},
\end{equation}
where the inequalities of the right side are considered in $p\A1 p$ and
$(1-p)\A1 (1-p)$, respectively. 
Moreover, the fibration $\varphi $ splits by means of
the $C^\infty$ global cross section 
\begin{equation}\label{4.8}
s :\qa \to  \pa \times G^+ \ , \quad \hbox { given by } \quad 
s (q) =(\phi (q) ,  |2q-1 |) ,
\end{equation}
for $q \in \qa$, where $|z|=(z^*z)^{1/2}$.
\end{pro}
\smallskip
\noi Proof. \rm  Let us first verify (\ref{4.7}). Fix $q\in \qa$. 
The only possible first coordinate of every pair in $\varphi \inv (q)$ must be 
$p = \phi (q)$, since it is the unique projection in $\pa $ 
with the same range as $q$. 

Given $ a = 
\left( \begin{array}{cc} a\sb 1 &a\sb 2  \\ 
a\sb 2^*& a\sb 3 \end{array} \right) \in G^+$, 
we know by (\ref{3.5}) that 
$\varphi (p,a) = q$ if and only if $  a\sb 2 = a\sb 1 (q-p) = a_1 x$. Then  
$a = \left( \begin{array}{cc} a\sb 1 &a\sb 1 x  \\ 
x^* a\sb 1& a\sb 3 \end{array} \right) ) $.
The inequalities $ x^* a\sb 1 x < a\sb 3$ in $ (1-p)\A1 (1-p)$ 
and $a_1 >0 $ in $p\A1 p$ are  
easily seen to be equivalent to the fact that 
$\left( \begin{array}{cc} a\sb 1 &a\sb 1 x  \\ x^* a\sb 1& a\sb 3 \end{array} \right) \in G^+$. 
This shows (\ref{4.7}).

Denote by $\eps = 2q-1$. It is clear that $\eps ^2 = 1$, i.e. $\eps $ is 
a symmetry. Consider its polar decomposition $\eps = \rho \la $, where
$\la = |\eps | \in G^+ $ and $\rho $ is a unitary element of $\A1$. 
>From the uniqueness of the polar decomposition it follows
 that $\rho = \rho ^* = \rho \inv $,  i.e. $\rho$ is a unitary selfadjoint 
symmetry. Then, since $q = \frac{\eps +1}{2}$, 
$$
\la \inv q^* \la = \la \inv \left( \frac{\eps ^* +1}{2} \right) \la =
\la \inv \left( \frac{\la \rho +1}{2} \right) \la =
 \frac{ \rho \la +1}{2} = \frac{\eps +1}{2} = q .
$$
Therefore $q \in {\cal S}\sb \la (\A1)$ and $\vfi (p, \la ) = 
\vfi (s (q) ) = q$ proving that $s$ is a cross section of $\vfi$ \QED

\bigskip
\begin{num}\label{cin} \rm
The space $\pa$ is the selfadjoint part of the space $\qa$.  
But each $a\in G^+$ induces the star $\*a$ and therefore 
another submanifold of $\qa$ of $a$-selfadjoint idempotents.
Let $a \in G^+$ and  denote the $a$-selfadjoint part  of $\qa$ by 
\begin{equation}\label{4.9}
\paa \ = \paa (\A1 ) \ = \ \{ \ q\in \qa \ : \ q^{\*a} = q \ \} . 
\end{equation}
We are going to relate the manifolds $\pa$ and $\paa$.
There is an obvious way of mapping $\pa $ onto $\paa$, namely 
$p \mapsto a^{-1/2}pa^{1/2} $. Its tangent map is the restriction
of the isometric isomorphism $X \mapsto a^{-1/2}Xa^{1/2} $ from 
$\aut $ onto $\auta$ mentioned in section 2. We shall study some 
less obvious maps between $\pa$ and $\paa$. 

For a fixed $a \in G^+$, consider the map
\begin{equation}\label{4.10}
\fa : \pa \to \paa \quad \hbox { given by } \quad \fa (p ) = \varphi (p, a)
, \ p\in \pa . 
\end{equation}
Then $\fa$ is a diffeomorphism between the submanifolds $\pa$ and $\paa$
of $\qa$ and $\fa \inv $ is just the map $\phi $ of (\ref{4.5}) restricted 
to $\paa$. The problem which naturally arises is the study of the
tangent map of $\fa $ in order to compare different $\paa$, $a \in G^+$. 

\medskip
The tangent space $(T \paa )\sb q $ for $q \in \paa$ can be described as in 
(\ref{TP}), 
\begin{equation}\label{4.11}
(T \paa )\sb q = \{ Y =
\left( \begin{array}{cc} 0&y \\ y^{\*a} &0 \end{array} \right) \in \A1\} = 
T(\qa)\sb q \cap \auta , 
\end{equation}
where the matricial representations are in terms of $q$. Therefore any 
$Y \in (T \paa )\sb q $ is characteri\-zed by its 1,2 entry $y = q Y$ by the 
formula 
\begin{equation}\label{4.12}
Y = y + y^{\*a} = qY + Yq . 
\end{equation}
\medskip  
\end{num}
\begin{pro}\label{4.13}
Let $p\in\pa $, $a \in G^+$ and 
$X = \left( \begin{array}{cc} 0&x \\ x^* &0 \end{array} \right) \in (T \pa )\sb p $. 
Denote by  $q = \fa (p) \in \paa$.  Then, in terms of $p$,  
$$
q \ (T\fa )\sb p (X) = 
\left( \begin{array}{cc} 0& a\sb 1 \inv x (a\sb 3 - a\sb 2^* a\sb 1 \inv a\sb 2)  \\ 
0 &0 \end{array} \right) = y .
$$
Therefore $(T\fa )\sb p (X) = y + y^{\*a}$ and $\|(T\fa )\sb p (X)\|\sb a = 
\|y \|\sb a$. 
\end{pro}
\smallskip
\noi Proof. \rm  We have the formula of Proposition \ref{3.3},
$$
\fa (p)= \fp(a)  = p E\sb p(a)\inv a = p\ (\ pap + (1-p)a(1-p)\ )\inv \ a .
$$
By the standard method of taking a smooth curve $\gamma$ in $\pa$ such 
that $\gamma (0) = p$ and $\dot \gamma (0 ) = X$, one gets 
$$
(T\fa )\sb p (X) = \left[ X- p E\sb p(a)\inv (Xap + paX - Xa(1-p)
-(1-p)aX ) \right]
E\sb p(a)\inv a \ . $$
Since $p$ and $E\sb p(a)$ commute, $pE_p(a)\inv (1-p)aX = 0$. In matrix form
in terms of $p$, by direct computation it follows that 
$$
\begin{array}{ll}
(T\fa )\sb p (X) = & \\
& \\
 =\left[ \ \left( \begin{array}{cc} 0&x \\ x^* &0 \end{array} \right) -
 \left( \begin{array}{cc} a\sb 1 \inv &0 \\ 0 &0 \end{array} \right) 
\left( \begin{array}{cc} xa\sb 2^* +a\sb 2 x^* &  a\sb 1x - xa\sb 3 \\
0 & 0 \end{array} \right) \ \right] 
\left( \begin{array}{cc} 1&a\sb 1 \inv a\sb 2  \\ 
a\sb 3 \inv a\sb 2 ^* &1 \end{array} \right) & \\
& \\
 = \left( \begin{array}{cc} 
-a\sb 1\inv x a\sb 2^* -a\sb 1\inv a\sb 2 x^*  & a\sb 1 \inv xa\sb 3  \\ 
x^*  &  0  \end{array} \right)  
\left( \begin{array}{cc} 1&a\sb 1 \inv a\sb 2  \\ 
a\sb 3 \inv a\sb 2 ^* &1 \end{array} \right) & \\
& \\
= \left( \begin{array}{ccc} -a\sb 1\inv a\sb 2 x^* & & 
a\sb 1 \inv (xa\sb 3 - xa\sb 2^*a\sb 1\inv a\sb 2  -a\sb 2  
x^*a\sb 1\inv a\sb 2) 
\\ x^* & & x^*a\sb 1\inv a\sb 2  \end{array} \right). & 
\end{array} 
$$
Multiplying by $q = \varphi (p,a)$, by (\ref{3.5}), one 
obtains
$$
\begin{array}{ll}
q \ (T\fa )\sb p (X)  = & \\
&\\
= \left( \begin{array}{cc} 1&a\sb 1 \inv a\sb 2  \\ 0 &0 \end{array} \right)
\left( \begin{array}{ccc} -a\sb 1\inv a\sb 2 x^* & & 
a\sb 1 \inv (xa\sb 3 - xa\sb 2^*a\sb 1\inv a\sb 2  -a\sb 2 x^*a\sb 1\inv a\sb 2) 
\\ x^* & & x^*a\sb 1\inv a\sb 2  \end{array} \right) & \\ 
& \\  
= \left( \begin{array}{cc} 0& a\sb 1 \inv 
x (a\sb 3 - a\sb 2^* a\sb 1 \inv a\sb 2)  \\ 
0 &0 \end{array} \right) = y , & 
\end{array}
$$
as desired. The fact that $\|y\|\sb a = \|Y\|\sb a$ is clear by regarding 
them as elements of  $(\A1 , \*a ) $ and using (\ref{4.11}) .\quad  

\bigskip

\section{The polar decomposition.}
\medskip
In this section it is convenient to identify $\qa$ with the set of 
symmetries (or reflections) $\{ \eps \in \A1 : \eps ^2 = 1 \}$ and 
$\pa $ with the set of selfadjoint symmetries
$\{ \rho \in \A1 : \rho = \rho ^* = \rho \inv \}$ by means of the 
affine map $x \mapsto 2x-1$. 

Recall that every invertible element $c$ of a unital \csta \ admits 
polar decompositions $c = \rho _1 \la _1 = \la _2 \rho _2 $, with
$\la_1, \la_2 \in G^+$ and $\rho_1, \rho_2 \in \U $. 
Moreover, 
$$
\la_1 = |c|, \ \la_2 = |c^*| \quad \hbox{  and  } \quad 
\rho_1 = \rho_2 = |c^*|\inv c = c|c|\inv . 
$$
In particular, if $\eps $ is a symmetry, its polar decompositions are 
$\eps = |\eps^*| \rho = \rho |\eps | $ and
\begin{equation}\label{rho}
\rho = \rho ^* = \rho \inv \in \pa .
\end{equation}
This remark defines the retraction
\begin{equation}\label{5.1}
\proj : \qa \to \pa , \quad \hbox{ by } \ \proj (\eps ) = \rho . 
\end{equation}
The map $\proj $ has been studied from a differential geometric
viewpoint in \cite{[CPR1]}. If $\eps \in \qa$, it is easy to show that 
$|\eps ^*| = |\eps |\inv$ and $|\eps ^*|^{1/2} \rho = 
\rho |\eps ^*|^{-1/2}$ (see \cite{[CPR1]}). This section is devoted to 
study,  for each $a\in G^+$, the restriction 
\begin{equation}\label{5.4}
\proj \sb a = \proj |\sb {\paa} : \paa \to \pa . 
\end{equation}
Observe that, with the identification mentioned above, 
$\paa = \qa \cap \auta = \qa \cap {\cal U} _a $, where
${\cal U} _a = \{u \in G : u\inv = u^{\*a} \}$ is the group of $a$-unitary
elements of $\A1$.
\medskip
\begin{pro}\label{5.5}
For every $a \in G^+$ the map 
$\proj \sb a : \paa \to \pa$  of (\ref{5.4}) is a diffeomor\-phism. 
\end{pro}
\smallskip
\noi Proof. \rm  By the remarks above, for every $\eps \in \qa$  
\begin{equation}\label{5.6}
\proj \sb a (\eps) = \rho = |\eps | \eps 
\end{equation}
which is clearly a $C^\infty$ map. 

Set $b = a ^{1/2}$ and consider, for a fixed $\rho \in \pa$,  
the polar decomposition of $b\rho b $ given by $ b\rho b = w |b\rho b |$, 
with $w \in \U$. Since $ b\rho b $ is invertible and selfadjoint 
by (\ref{rho}), it is easy to prove 
(see \cite{[CPR3]}) that
$$
w = w^* = w\inv \in \pa \ , \quad w\ b\rho b = b\rho b \ w 
\quad \hbox { and } \quad 
w \ b\rho b = |b\rho b | \in G^+ . 
$$
Let $\eps = b\inv w b$. It is clear by its construction that 
$\eps \in \paa $. Also $\eps \rho = \la >0 $, since  
$$
b \eps \rho b = w \ b\rho b = |b\rho b | \in G^+ . 
$$
Therefore the polar decomposition of $\eps $ must be $\eps = \la \rho $. 
So $\la = |\eps ^*|$ and $\proj \sb a (\eps) = \rho$. Therefore 
\begin{equation}\label{5.7}
\proj \sb a \inv (\rho ) = a^{-1/2} \left( a^{1/2}\rho a^{1/2} 
|a^{1/2}\rho a^{1/2}|\inv \right) a^{1/2} , 
\end{equation}
which is also a $C^\infty $ map, showing that $\proj \sb a$ is a 
diffeomorphism. \quad \QED
\medskip
\begin{rem}\label{5.8}\rm
The fibres of the retraction $\proj $ over each $p \in \pa $  are 
in some  sense, ``orthogonal'' to $\pa$. 
In order to explain this remark, consider the algebra $\A1 = M_n (\zC )$ 
of all $n \times n$ matrices with complex entries. Then $M_n (\zC )$ has 
a natural scalar product given by $<X,Y> = tr (Y^*X)$. It is easy to
prove that for every $p \in \pa $, $(T\pa )_p $ is orthogonal to 
 $(T\proj \inv (\pa ))_p$. The same result holds in every
\csta \ with a trace $\tau$.
Then the map $a \mapsto \proj \sb a (\rho) \inv $ of 
(\ref{5.7}) can be considered as the ``normal'' movement which produces 
$a$-selfadjoint projections for every $a \in G^+$.

On the other hand, the map $\fp $ of 
(\ref{3.5}), which was studied also in \cite{[PW1]}, gives another way to get 
$a$-selfadjoint
projections for every $a \in G^+$. In terms of the geometry of $\qa$ this
way is, in the sense above, an oblique movement. 
A related  movement is to take for each $a \in G^+$, an $a$-selfadjoint 
projection $q'$ with $\ker q' = \ker p $.

\end{rem}
\bigskip
Combining, for a fixed $a \in G^+$ the maps $\fa $ of (\ref{4.10}) and 
$\proj \sb a$ of (\ref{5.4}), one obtains 
a $C^\infty$ movement of the space $\pa$. The following proposition 
describes explicitly this movement.
\medskip
\begin{pro}\label{5.9}
Let $a \in G^+$. Then the map 
$\proj \sb a \circ \fa : \pa \to \pa$  is a diffeomor\-phism 
of $\pa $. For $p \in 
\pa $, let $\fa (p) = q = p+x$ and $\eps = 2q-1 $. In terms of $p$,  
$x = a\sb 1 \inv a\sb 2 $ if $a = \left( \begin{array}{cc} a\sb 1 &a\sb 2 \\ 
a\sb 2^*& a\sb 3 \end{array} \right) $ and 
\begin{equation}\label{5.10}
\begin{array}{rl}
 \proj \sb a \circ \fa (p)                   &  =  \ 
\left( \begin{array}{cc} 1+xx^* &0  \\0& 1+x^*x \end{array} \right) ^{-1/2}
\left( \begin{array}{cc} 1 &x  \\x^*& -1 \end{array} \right) \\
&  \\
& = \ \left[ qq^* + (1-q)^*(1-q) \right] ^{-1/2} (q+q^*-1). 
\\
\end{array}
\end{equation}
\end{pro}
\smallskip
\noi Proof. \rm  In matrix form, $\eps = 2\fa (p)-1 = 
\left( \begin{array}{cc} 1 &2x  \\0& -1 \end{array} \right)$ so that 
$$
\eps ^* \eps = \left( \begin{array}{cc} 1 &0  \\2x^*& -1 \end{array} \right)
\left( \begin{array}{cc} 1 &2x  \\0& -1 \end{array} \right) = 
\left( \begin{array}{cc} 1 &2x  \\2x^*& 4x^*x+ 1 \end{array} \right) = | \eps |^2 .
$$
On the other hand, by (\ref{4.8}),  $q \in P_{| \eps |}(\A1 )$. 
Therefore, by (\ref{3.5}), \\
$| \eps | = 
\left( \begin{array}{cc} b &bx  \\x^*b& c \end{array} \right) $ 
with $b, \ c $ positive. Straightforward computations show that 
$$
b = (1+xx^*)^{-1/2}\quad \hbox { and } \quad 
c^2 = 4x^*x +1 - x^* (1+xx^*)\inv x . 
$$
Since $x^*(1+xx^*) = (1+x^*x)x^*$, 
$$
\begin{array}{rl} 
c^2 & = 4x^*x+1 - (1+x^*x)\inv x^*x \\
& \\    
& = 4x^*x + (1+x^*x) \inv       \\
&\\    
& = (1+x^*x) \inv (4(x^*x)^2 + 4x^*x +1 ) \\
&\\    
& = (1+x^*x) \inv (2x^*x +1 )^2 .
\end{array}
$$
Then $c =  (1+x^*x)^{-1/2} (2x^*x+1)$ and
\begin{equation}\label{modeps}
| \eps | = \left( \begin{array}{cc} 1+xx^* &0  \\0& 1+x^*x \end{array} \right) ^{-1/2}
\left( \begin{array}{cc} 1 &x  \\x^*& 2x^*x+1  \end{array} \right) .
\end{equation}
Now the two formulas of (\ref{5.10}) 
follow by easy matrix computations.  \quad \QED
\begin{rem}\rm
It is interesting to observe that the factor $q+q^*-1$ of (\ref{5.10})
has been characterized by Buckholtz \cite{[Bu]} as the inverse of
$P_{R(q)}-P_{\ker q}$. 
\end{rem}

\bigskip
A natural question about these movements is the following: 
for $p \in \pa$, how far can 
$\proj \sb a \circ \fa (p) $ be from $p$. In order to answer this question
we consider the orbit 
\begin{equation}\label{5.11}
\op := \{ r \in \pa : \proj \sb a \circ \fa (p) = r \ \hbox{ for some } \ 
a \in G^+ \} . 
\end{equation}
The next result is a metric characterization
of $\op$  based on some results about the ``unit disk'' of the 
projective space of $\A1$ defined by $p$ (see \cite{[ACS]}). 
\medskip
\begin{pro}\label{5.12}
Let $p \in \pa$. Then 
$$
\op = \{ r \in \pa : \|r-p\| < \frac{\sqrt {2}}{2} \} 
$$
\end{pro}
\smallskip
\noi Proof. \rm Fix $a \in G^+$. Let $q = \fa (p)$,  $\eps = 2q-1$
and $r = \proj \sb a \circ \fa (p)$. By (\ref{4.8}), $r$ is also 
obtained  if we replace
$a$ by $|\eps |$, since $\fa (p) = q = \vfi \sb { |\eps |} (p) $ and 
$r = \proj \sb a(q) = \proj (q) = \proj \sb {|\eps |}(q)$. 
Note that $|\eps |$ is positive and $\rho$-unitary, i.e. unitary 
for the signed inner product $<,>\sb \rho$ given by  
$\rho = |\eps | \eps = 2r-1 $. Indeed, by (\ref{rho}),  
$|\eps |^{\#\sb {\rho}} = \rho \inv |\eps | \rho = |\eps |\inv $.

Since $\eps = \rho |\eps | = |\eps |^{-1/2} \rho |\eps |^{1/2} $, also
$q = |\eps |^{-1/2} r |\eps |^{1/2}$. In \cite{[ACS]} it is shown that 
the square root of a $\rho$-unitary is also $\rho$-unitary. Then
$|\eps |^{-1/2}$ is $\rho$-unitary. In \cite{[ACS]} it is also shown that
$$
\|r - P\sb {R(\la r \la \inv )} \| < \frac{\sqrt {2}}{2}
$$
for all positive $\rho$-unitary $\la$. Note that 
$p = P\sb {R(q)} $ and then it must be $\|p-r\|  < \frac{\sqrt {2}}{2}$.
In Proposition 6.13 of \cite{[ACS]}  it is shown 
that for all $r\in \pa $ such that 
$\|r-p\|  < \frac{\sqrt {2}}{2} $, there exists a positive $(2r-1)$-unitary
$\la$ such that $p = P\sb {R(\la r \la \inv )}$. In this case
$r = \proj \sb \la \circ \vfi \sb \la (p) \in \op$. \quad \QED

\bigskip
\section{New short geodesics.}
\medskip
Lengths of geodesics in $\pa$ have studied in \cite{[PR2]}, \cite{[Br]}, 
\cite{[Ph]} and \cite{[ACS]}.
It has been proved that if $p,r \in \pa $ and 
$\|p-r\| <1 $, then there exists a unique geodesic of $\pa$ 
joining them which has minimal length. 
On the other hand, the fibres $\proj \inv (\pa )$ are geodesically complete
and the geodesic joining $q_1, q_2 \in \proj \inv (\pa )$ is a shortest 
curve in $\qa$ \cite{[CPR2]}. This final section is devoted to show the existence 
of ``short oblique geodesics'', i.e. geodesics which are not contained 
neither in $\pa$ nor in the fibres. 

More precisely, the idea of the present section is to 
the use different stars $\*a$ for 
$a \in G^+$ in order to find short curves between  pairs of 
non  selfadjoint idempotents of $\A1$. Basically we want to characterize 
those pairs $q, r \in \qa$ such that there exist $ a \in G^+$ with
$q,r \in \paa$. If $q$ and $r$ remain close in $\paa$, they can be
joined by a short curve in the space $\paa$.

The first problem is that the positive $a$ need not be unique. This
can be fixed up in the following manner:
\medskip
\begin{lem}\label{unigeo} 
Suppose that $a \in G^+$ and $p,r \in \pa \cap \paa$. Then
$\|p-r\| = \|p-r\|\sb a$ and, if $\|p-r\|<1 $, the short geodesics 
which join them in $\pa $ and $\paa $ are the same and have 
the same length. 
\end{lem}
\dem 
Note that $\pa \cap \paa$ is the space of projections commuting 
with $a$. Let $\B1 = \{a\}' \cap \A1 $, the relative commutant of $a$ in 
$\A1$. Since $a = a^*$, $\B1$ is a \csta . Moreover, $\pa \cap \paa = \pb $.
Now, since $\|p-r\|<1$, $p$ and $q$ can be joined by the unique short geodesic $\gamma$ along $\pb$ (see \cite{[PR2]} or \cite{[ACS]}) and $\gamma$ is also a 
geodesic both for $\pa $ and $\paa $. The length of $\gamma$ is
computed in the three algebras in terms of the $norm$ of the 
corresponding tangent vector $X$. But since $X\in \B1$, its norm is the
same with the two scalar products involved.   \quad \QED

\medskip
We shall give a characterization of pairs of close idempotents 
$p, q \in \qa$ such 
that $p,q \in \paa $ for some $a \in G^+$. The characterization will be 
done in terms of a tangent vector $X \in T(\qa )_p $ such that 
$q = e^X p e^{-X}$. First we give a slight improvement of the way 
to obtain such $X$ which appears in 2) of \cite{[PR2]}: 
\medskip
\begin{pro}\label{el X}
Let $p \in \pa $ and $q\in \qa$ with $\| p-q\| < 1$. Let $\eps = 2q-1$, 
$\rho = 2p-1$, 
$$
v_1 = \frac{\eps \rho +1}{2} = qp +(1-q)(1-p) \quad \hbox {  and  }
\quad v_2 = \frac{\rho \eps +1}{2} = pq +(1-p)(1-q) .
$$ 
Then $\|v_1 -1\| = \|v_2 -1\|=\|p-q\| < 1 $ and
\begin{equation}\label{log}
X = (Id - E_p)(\log v_1 ) =  \frac12 \ ( \log v_1 - \log v_2 ) 
\end{equation}
verifies that $X \in T(\qa)_p $ (i.e. $pXp = (1-p)X(1-p) =0 $) and 
$q = e^X p e^{-X}$. 
\end{pro}
\dem
Note that $\rho = \rho ^* = \rho \inv  \in \pa$. Then 
$$
\|v_1 -1\| = \| \frac{\eps \rho - 1 }{2} \| = 
\frac12 \| (\eps - \rho )\rho \|=\| q-p\| < 1 , 
$$
and similarly for $v_2$. Let $X_i = \log v_i $ for $i = 1, 2$. Since 
$v_1 \rho = \rho v_2$, and each $X_i$ is obtained as a power 
series in $v_i$, 
we obtain also that $X_1 \rho = \rho X_2$. Then, if $X= \frac12 (X_1 - X_2)$, 
we have that $X\rho = -\rho X$, and then $X\in T(\qa )_p$. 

Note also that $v_1 \rho = \eps v_1$ and $\|v_i -1\|<1$ for $i = 1,2$. So
$v_1 p v_1 \inv = q$.  Easy calculations show that  $v_1$ and $v_2 $
commute. As before this implies that $X_1$ and $X_2$ commute. 
Then
$$
w = v_1 v_2 = v_2 v_1 = e^{X_1 +X_2} = ( \frac{\eps +\rho}{2} ) ^2
$$ 
commutes with $v_1 , \ v_2 ,  \ \rho , \ \eps ,\ p$ and $q$. Denote by 
$$
w^{-\frac12} = e^{- \frac{X_1 +X_2}{2}} . 
$$
Since $(X_1 +X_2 )\rho = \rho (X_1 +X_2 )$, $w^{-1/2}$ commutes with 
$\rho$. Note that 
\begin{equation}\label{elX}
X = \frac{X_1-X_2}{2} = X_1 -\frac{X_1 + X_2}{2}.
\end{equation}
This implies that $e^X= e^{X_1} w^{-1/2} = v_1 w^{-1/2}$ and 
therefore 
$$
e^Xpe^{-X} = v_1 p v_1 \inv = q .
$$   

Finally, since $X$ has zeros in its diagonal 
and $\frac{X_1 + X_2}{2}$ is diagonal
in terms of $p$, we deduce from (\ref{elX}) 
that $X = (Id - E_p)(X_1)$ and the proof 
is complete. \quad \QED

\medskip
\begin{pro}\label{6.1}
Let $p\in \pa$ and $q \in \qa$ such that 
$\|p-q\|<1$. Let $X = 
\left( \begin{array}{cc} 0 &x  \\y& 0 \end{array} \right) \in T(\pa )\sb p$ 
as in (\ref{log}), such that
$e^X p e^{-X} = q$. Then the following are 
equivalent: 
\ben
\item  There exists $a \in G^+$ such that $p, q \in \paa$.
\item  There exists $a \in G^+$ such that $pa=ap$ and $X^{\*a} = -X$.
\item  There exist $b \in G^+(p\A1 p) $ and $c \in G^+((1-p)\A1 (1-p) )$ 
such that 
$$y = - cx^*b .$$
\een
\end{pro}
\smallskip
\noi Proof. \rm Condition 2. can be written as
$$
a = \left( \begin{array}{cc} a\sb 1 &0  \\0& a\sb 2 \end{array} \right) \quad
\hbox{ and } \quad a\inv X ^* a = -X .
$$
In matrix form
$$
\begin{array}{rl}
a\inv X ^* a & = 
 \left( \begin{array}{cc} a\sb 1\inv  &0  \\0& a\sb 2\inv  \end{array} \right)
 \left( \begin{array}{cc} 0 &y^*  \\x^*&0  \end{array} \right)
 \left( \begin{array}{cc} a\sb 1 &0  \\0& a\sb 2 \end{array} \right) \\
&   \\
&  = 
\left( \begin{array}{cc} 0 &a\sb 1\inv y^*a\sb 2  
\\a\sb 2\inv x^*a\sb 1  & 0 \end{array} \right) \\
&  \\
& =  \left( \begin{array}{cc} 0 &-x  \\-y&0  \end{array} \right) , 
\end{array}
$$
which clearly is equivalent to condition  3.

Condition 1. holds if $X^{\*a} = -X$, since in that case 
$e^X$ is $a$-unitary and then $q \in \paa$. In order to prove the converse,
we consider $\*a $ instead of $*$ and so 
condition 1. means that $p,q \in \pa$.
Then, with the notations of (\ref{el X}), 
we have that $v_2 = v_1 ^*$ and  
$$
X_2 = \log v_2 =\log v_1^* = X_1^* \quad \Rightarrow  \quad 
X^* = \frac{X_2 - X_1}{2} = -X , 
$$ 
showing 2. \quad \QED
\medskip
\begin{rem} \rm \

\smallskip
Let us call a direction (i.e. tangent vector) in $\qa$  ``good"
if it is the direction of a short geodesic. Proposition \ref{6.1} 
provides a way to obtain good directions. Other good directions occur in the 
spaces $p \A1 (1-p)$ and $(1-p) \A1 p$, 
determined by the affine spaces 
of projections with the same range ($\qpa $) 
or the same kernel as $p$, where the straight lines can be 
considered as short geodesics.

Still another good directions can be found
looking at pairs $p, q \in \qa$ such that, for some $a\in G^+$,  
$\proj^a (q) = p$, where $\proj^a$ means the retraction of (\ref{5.1}), 
considering in $\A1$ the star $\*a$. These pairs can be characterized
in a very similar way as Proposition \ref{6.1}. In fact, in condition 3 
(with the same notations), $y = -bx^*c$ should be replaced by $y =bx^*c$.
These directions are indeed good because it is known \cite{[CPR2]} that
along the fibers of each $\proj ^a$  there are short geodesics 
that join any pair of elements (not only close pairs). 
\end{rem}

\vglue1truecm
\direcciones

\end{document}